\documentclass[11pt, a4paper]{amsart}
\usepackage{palatino}
\usepackage{amsmath,amssymb,amsfonts}
\usepackage[margin=3cm]{geometry}
\usepackage[breaklinks,bookmarksopen,bookmarksnumbered]{hyperref}

\theoremstyle{plain}
\newtheorem*{thm}{Theorem}
\newtheorem*{lem}{Lemma}

\newtheorem*{cor}{Corollary}
\theoremstyle{remark}
\newtheorem*{rem}{Remark}

\def\pr{\noindent\textit{Proof} : }

\def\rond{\kern 1pt{\scriptstyle\circ}\kern 1pt}

\newcommand\Aut{\operatorname{Aut}}

\def\Z{\mathbb{Z}}

\def\C{\mathbb{C}}
\def\P{\mathbb{P}}

\def\F{\mathbb{F}}

\def\Sy{\mathfrak{S}}
\def\A{{\mathfrak{A}}}
\def\H{\mathrm{H}}

\def\iso{\vbox{\hbox to .8cm{\hfill{$\scriptstyle\sim$}\hfill}
\nointerlineskip\hbox to .8cm{{\hfill$\longrightarrow $\hfill}} }}

\catcode`\@=11\def\eqalign#1{\null\,\vcenter{\openup\jot\m@th\ialign{
\strut\hfil$\displaystyle{##}$&$\displaystyle{{}##}$
&&\quad\strut$\displaystyle{##}$&$\displaystyle{{}##}$
\crcr#1\crcr}}\,}
\catcode`\@=12

\begin{document}
\title[Non-rationality of the symmetric sextic Fano threefold]{Non-rationality of the symmetric sextic Fano threefold}
\author[Arnaud Beauville]{Arnaud Beauville}
\address{Laboratoire J.-A. Dieudonn\'e\\
UMR 6621 du CNRS\\
Universit\'e de Nice\\
Parc Valrose\\
F-06108 Nice cedex 2, France}
\email{arnaud.beauville@unice.fr}
 
\date{\today}

\maketitle

\section*{Introduction}

The symmetric sextic Fano threefold is the subvariety  $X$ of $\P^6$ defined by the equations:
\[\sum X_i=\sum X_i^2=\sum X_i^3=0\ .\]
It is a smooth  complete intersection of a quadric and a cubic in $\P^5$, with an action of $\Sy_7$.
We will prove that it is not rational.
\par Any smooth complete intersection of a quadric and a cubic in $\P^5$ is unirational  \cite{E}.
 It is known that a \emph{general} such intersection is not rational: this is proved in \cite{B} (thm. 5.6) using the intermediate Jacobian, and in \cite{Pu} using the group of birational automorphisms. But neither of these methods allows to prove the non-rationality of any particular such threefold.

\par Our motivation comes from the recent paper of Prokhorov \cite{P}, which classifies the simple finite subgroups of the Cremona group $\mathrm{Cr}_3=\mathrm{Bir}(\P^3)$. In view of this work our result implies that the alternating group $\A_7$ admits
 only one embedding into $\mathrm{Cr}_3$ up to conjugacy.

\par To prove our result we use the Clemens-Griffiths criterion (\cite{CG}, Cor. 3.26):  if $X$ is rational, its intermediate Jacobian $JX$ is the Jacobian of a curve, or a product of such Jacobians. The presence of the automorphism group $\Sy_7$,
together with the celebrated bound $\ \#\Aut(C)\leq 84(g-1)$ for a curve $C$ of genus $g$,  immediately implies that $JX$ is not isomorphic to the Jacobian of a curve. To rule out products of Jacobians we need some more information, which is provided by
 a simple analysis of the representation of $\Sy_7$ on the tangent space $T_0(JX)$.

\medskip
\section*{The result}
\begin{thm}
The intermediate Jacobian $JX$ is not isomorphic to a Jacobian or a product of Jacobians. As a consequence, $X$ is not rational.
\end{thm}
That the second assertion follows from the first is the Clemens-Griffiths criterion mentioned in the introduction. Since the
Jacobians and their products form a closed subvariety of the moduli space of principally polarized abelian varieties, this gives an easy proof of the fact that a 
general intersection of a quadric and a cubic in $\P^5$ is not rational.

\par As mentioned in the introduction, the classification in \cite{P} together with the theorem   implies:

\begin{cor}
Up to conjugacy, there is only one embedding of $\mathfrak{A}_7$ into the Cremona group $\mathrm{Cr}_3$, given by an embedding $\mathfrak{A}_7\subset \mathrm{PGL}_4(\C)$.\qed
\end{cor}
(The embedding $\mathfrak{A}_7\subset \mathrm{PGL}_4(\C)$ is the composition of the standard representation $\mathfrak{A}_7\rightarrow 
\mathrm{SO}_6(\C)$ and the double covering $\mathrm{SO}_6(\C)\rightarrow  \mathrm{PGL}_4(\C)$.)

\medskip
\par  The intermediate Jacobian $JX$ has dimension 20. The group $\Sy_7$ acts on $JX$ and therefore on the tangent space $T_0(JX)$; we will first determine this action. 

\begin{lem}
As a $\Sy_7$-module $T_0(JX)$ is the sum of two irreducible representations, of dimensions $6$ and $14$.
\end{lem}
\pr Let $V$ be the standard (6-dimensional) representation of  $\Sy_7$, and $\P:=\P(V)$; we will view $X$ as a subvariety of $\P$, stable under $\Sy_7$.
\par Recall that $T_0(JX)$ is  $\H^1(X,\Omega ^2_X)$, and that  the exterior product $\Omega ^1_X\otimes \Omega ^2_X\rightarrow K_X$ induces a canonical isomorphism $\Omega ^2_X\iso T_X(-1)$. 
The exact sequence
\[0\rightarrow T_X \longrightarrow T_{\P}{}^{}_{|X}\longrightarrow \mathcal{O}_{X}(2)\oplus \mathcal{O}_{X}(3)
\rightarrow 0\]
twisted by $\mathcal{O}_X(-1)$, gives a cohomology exact sequence
\[ \hss 0\rightarrow \H^0(X, T_{\P}(-1)^{}_{|X}) \rightarrow \H^0(X,\mathcal{O}_{X}(1))\oplus  \H^0(X,\mathcal{O}_{X}(2)) \rightarrow  \H^1(X,T_X(-1)) \\
\longrightarrow   \H^1(X, T_{\P}(-1)^{}_{|X})\ .\hss\]

From the Euler exact sequence $\ 0\rightarrow \mathcal{O}_X\rightarrow  \mathcal{O}_X(1)\otimes_{\C}V \rightarrow T_{\P}{}^{}_{|X}\rightarrow 0\ $ we deduce  $ \H^1(X, T_{\P}(-1)^{}_{|X})=0$  and an isomorphism $V\iso  \H^0(X, T_{\P}(-1)^{}_{|X})$.
Thus we find an exact sequence
\[0\rightarrow V \longrightarrow   \H^0(X,\mathcal{O}_{X}(1))\oplus \H^0(X,\mathcal{O}_{X}(2))  \longrightarrow T_0(JX)\rightarrow 0\ , \]
which is equivariant with respect to the action of $\Sy_7$. As $\Sy_7$-modules $ \H^0(X,\mathcal{O}_{X}(1))$ is isomorphic to $V$, and $ \H^0(X,\mathcal{O}_{X}(2))$ to  $\mathsf{S}^2V/\C.q$, where $q$ corresponds to the quadric containing $X$. On the other hand we have 
$\mathsf{S}^2V=\C\oplus V\oplus V_{(5,2)}$, where $V_{(5,2)}$ is the irreducible representation of $\Sy_7$ corresponding to the partition $(5,2)$ of $7$ (\cite{FH}, exercise 4.19). Thus we get $T_0(JX)\cong V\oplus V_{(5,2)}$.\qed

\bigskip
\noindent\emph{Proof of the theorem} :
We first observe that $\Sy_7$ cannot act faithfully on the Jacobian $JC$ of a curve of genus $g\leq 20$. Indeed by the Torelli theorem, the map $\Aut(C)\rightarrow \Aut(JC)$ is injective and its image has index 1 (if $C$ is hyperelliptic) or  2 otherwise. Thus we find $\# \Aut(C)\geq \frac{1}{2}7!=2520$. On the other hand we have $\# \Aut(C)\leq 84(g-1)\leq 1596$, a contradiction.\label{jac}

\smallskip
 Let $A$ be a principally polarized abelian variety. Recall that $A$ can be written \emph{in a unique way} as a product
$A_1^{a_1}\times \ldots \times A_p^{a_p}$, where  $A_1,\ldots ,A_p$ are indecomposable, non isomorphic principally polarized abelian varieties.
(this decomposition corresponds to the decomposition of the Theta divisor into irreducible components, see \cite{CG}, Cor. 3.23). Therefore we have 
\[ \Aut(A)\cong \Aut(A_1^{a_1})\times \ldots \times \Aut(A_1^{a_1})\quad \hbox{and}\quad \Aut(A_i^{a_i})\cong \Aut(A_i)^{a_i} \rtimes \Sy_{a_i}\ . \]
Let $G\subset \Aut(A)$; the product decomposition of $A$ induces a decomposition of $G$-modules
\[T_0(A)=T_0(A_1^{a_1})\oplus\ldots \oplus T_0(A_p^{a_p})\ .\]
Fix an integer $i$; the group $G$ permutes the factors of $A_i^{a_i}$. 
If $O_1,\ldots ,O_k$ are the orbits of $G$ in this action, we have a further decomposition of $G$-modules
\[T_0(A_i^{a_i})=T_0(A_i^{O_1})\oplus\ldots \oplus T_0(A_i^{O_k})\ .\]

\smallskip
 Assume that   $JX$ is isomorphic to a product of Jacobians $A_1^{a_1}\times \ldots \times A_p^{a_p}$, with $A_i\not\cong A_j$ for $i\neq j$. If $\Sy_7$ acts on a set with $\leq 20$ elements, its orbits have order $1$ or $7$ (\cite{DM}, thm. 5.2.B). An orbit with one element  means that $\Sy_7$ acts  on the Jacobian $A_i$; by the lemma  this action is faithful,
 contradicting the beginning of the proof. Thus each $a_i$ must be divisible by 7, which is impossible since $\sum a_i\dim(A_i)=20$.\qed

\bigskip

\begin{rem}
The same kind of argument gives a simple proof that the Klein cubic threefold, defined by $\sum_{i\in \Z/5}X_i^2X_{i+1}=0\ $ in $\P^4$, is not rational (and by the same token that the general cubic threefold is not rational). The automorphism group of the Klein cubic is $\mathrm{SL}_2(\F_{11})$, of order 660, while its intermediate Jacobian has dimension 5. It is easily seen as above that a 5-dimensional principally polarized abelian variety with an action of $\mathrm{SL}_2(\F_{11})$ cannot be a Jacobian or a product of Jacobians (see also \cite{Z} for a somewhat analogous, though more sophisticated, proof).
\end{rem}

\end{document}